\documentclass[11pt,psfig]{article}
\usepackage{amsmath,amssymb,epsfig}

\newtheorem{theorem}{Theorem}[section]
\newtheorem{definition}[theorem]{Definition}

\newtheorem{lemma}[theorem]{Lemma}

\newtheorem{remark}[theorem]{Remark}

\def\sF{{\cal F}}

\def\P{\mathbb{P}}
\def\Q{\mathbb{Q}}

\def\sF{{\cal F}}

\def\P{\mathbb{P}}

\makeatletter % so @ will be a letter
\@addtoreset{equation}{section}

\makeatother % so @ will not be a letter
\begin{document}
\begin{center}
{\Huge\bf On Robust Utility Maximization}
\mbox{}\\
\vspace{2cm}
Traian A. Pirvu\\
Department of Mathematics\\
The University of British Columbia\\
Vancouver, BC, V6T1Z2\\
tpirvu@math.ubc.ca\\
\vspace{1cm}
Ulrich G. Haussmann\footnote{Work
supported by NSERC under research grant 88051 and NCE grant 30354 (MITACS).}\\
Department of Mathematics\\
The University of British Columbia\\
Vancouver, BC, V6T1Z2\\
uhaus@math.ubc.ca\\

\mbox{}\\

\today
\end{center}

\noindent {\bf Abstract.}
This paper studies the problem of optimal investment in incomplete markets, robust with respect to stopping times.
We work on a Brownian motion framework and the stopping times are adapted to the Brownian filtration. Robustness
can only be achieved for logartihmic utility, otherwise a cashflow  should be added to the investor's wealth. The cashflow can be decomposed into the sum of an increasing and a decreasing process. The last one can be viewed as consumption. The first one is an insurance premium the agent has to pay.

\vspace{1cm}


\begin{quote}

\end{quote}

\begin{flushleft}
\end{flushleft}

\setcounter{equation}{0}
\section{Introduction}
Dynamic asset allocation has been an important field in modern finance. The ground-breaking paper in this literature is Merton \cite{Mer71}. He assumed a utility function of the power type (CRRA) and the market consisting
of a risk-free asset with constant rate of return and one or more stocks, each with constant mean rate of return and volatility. Merton was able to derive closed form solution for the stochastic control problem of maximizing utility of final wealth.

Karatzas et al. \cite{KarLehShr87}, Cox and Huang \cite{CoxHua89} establish the static martingale method in solving the optimal
investment and consumption problem in the complete market paradigm. In incomplete markets perfect risk transfer is not possible, thus people selected some \textit{martingale measures} according to some
risk criteria. He and Pearson \cite{HE} considered \textit{the minmax martingale measure} to transform the dynamic portfolio allocation problem into a static one. An auxiliary problem is analyzed, where the budget constraint is turned into a static constraint using only one martingale measure, i.e., the final wealth is feasible just under this martingale measure. \textit{The minmax martingale measure}
is defined to be the martingale measure for which the solution of the auxiliary problem, coincides with the solution of the original one. Kramkov and Schachermayer \cite{Kramkov} analysed the problem of maximizing utility from final wealth in a general semmimartingale model by means of duality.

This work combines the problem of utility maximization with the problem of hedging contingent claims in incomplete markets.   
In incomplete markets, when it comes to pricing and hedging a crucial issue is how to choose the appropriate martingale or risk-neutral measure. There are (among others) two main competing quadratic hedging approaches: local risk-minimization and mean-variance. They give rise to two martingale measures:
\textit{the minimal martingale measure} and \textit{the variance-optimal martingale measure} (see Heath et al. \cite{HEATH}).
These two martingale measures coincide if the stock price has independent increments (see Grandits \cite{GRA}). Moreover \textit{the minmax martingale measure} for quadratic utilities is also 
\textit{the variance-optimal martingale measure} (see Ex $5.3$ in \cite{BEL}). 

Our aim is to understand time consistency in an investor's optimal strategy. According to  Merton \cite{Mer71} an agent with CRRA preferences should invest a constant proportion of her wealth
in the risky assets. Is this strategy time consistent? The answer is NO, if due to some unforseen events such as death, or getting fired, the agent's investment horizon is some stopping time $\tau.$
Dynamic inconsistent behavior was first formalized analytically by Strotz \cite{STRO}. The problem arises if the investor at later dates is free to reconsider her policy. Assume
that her investment horizon is $T,$ but later on she learns it will be changed to some stopping time $\tau\leq T.$ Consequently she may choose to change her investment strategy.

This issue of horizon uncertainty  goes back to Levhari and Mirman \cite{LEV}. One can regard the stopping time horizon as a major event and by time change, it can be turned into a stochastic clock. Gol and Kallsen \cite{GOL} solve the problem of logarithmic utility maximization
in a general semimartingale framework and show that the optimal strategy is robust if one uses a stochastic clock. 
In a similar setup Bouchard and Pham \cite{BOU} extends duality techniques to characterize the optimal solution.
Blanchet-Scalliet et al \cite{BEL} considers a random horizon not necessarily a stopping time.
Zitkovi\'{c} \cite{ZIT} looks at the problem of utility maximization with a stochastic clock and an unbounded random endowment. Karatzas and Wang \cite{KWA} treats utility maximization problem
of mixed optimal stopping/control type.

  Choulli et all. \cite{CH} address the consistency problem in a semmimartingale setup when the preferences are CRRA. The core idea is that for a stopping time $\tau,$ the optimal investment of two agents (one can also think of two governments) over time intervals $[0,\tau]$ and $[\tau,T]$ is not the same as the investments of one agent over $[0,T]$ unless they are myopic. Strictly speaking this is saying that non-logarithmic utility maximization is not robust with respect to stopping times.
As pointed by  Choulli et al., the resolution is to add a correction term to investor's wealth.

In a different context Ekeland and Lazrak \cite{EKE} question time consistency of optimal consumption when discounting is not exponential. It turns out that the policy which is optimal from time
zero perspective will be still optimal at some later time $t$ only if the discounting is exponential. Thus if the investors at time zero cannot commit to the decision maker at later times $t>0,$
the optimal control approach will derive a policy which is impossible to implement. Instead they use a game-theoretic approach to derive an equilibrium strategy. 

Let us notice that the time consistency required by Choulli et al. is stronger because it involves all stopping times (not just the deterministic ones).

This work, in a Brownian framework, proposes a correction term process which needs to be added to an investor's wealth to achieve time consistency. The investor's risk preferences are more general than CRRA, thus it can be seen as an extension of Choulli et al.
The correction term is intimately related to the investor's coefficient of relative risk aversion and coefficient of prudence. The latter one was introduced by Kimball \cite{KIMB} as a measure of the sensitivity of choices to risk. In the complete market paradigm we establish the uniqueness of it. 

 Moreover we go one step farther and show how to finance it. If the coefficient of prudence is less than twice the coefficient
 of relative risk aversion the correction term is negative in which case the agent can just consume it. In general being a process of finite variation it can be decomposed
 as the sum of an increasing and a decreasing process. The increasing processes is an insurance premium the agent has to pay so that her policy is optimal if she would stop investing at any stopping time $\tau.$ The terminal value of it, which may be regarded as a contingent claim, can be represented by the Clark-Haussmann-Ocone formula and a hedging portfolio can be implemented. 

 Therefore to invest optimally and robustly with respect to stopping times, the agent should use some of her initial wealth to finance the hedging portfolio and consume the decreasing
 component and the difference between hedging portfolio and correction term. The same investment strategy can be implemented in incomplete markets. However since
 perfect risk transfer is no longer possible the agent would carry in her portfolio an intrinsic (unhedgeable) risk.  

Our main contribution to the field is finding a correction term process for risk preferences more general than CRRA and establish the uniqueness of it for the case of complete markets. It
would be interesting to characterize in incomplete market all the processes which added to investor's wealth lead to time consistency.

 The remainder of this paper is structured as follow. In Section $2$ we introduce the financial market model and section $3$ presents the objective. Sections $4$ and $5$ treats the complete and incomplete market case. Section $6$ derives the hedging portfolio. We conclude with an appendix containing some proofs.

\section{Model Description}
  
\subsection{Financial Market}
 We adopt a model for the financial market consisting of one bond and $d$ stocks. We work in discounted terms, that is the price
 of bond is constant and the stock price per share satisfy
 $$
dS_i(t)=S_i(t)\left[\alpha_i(t)\,dt
+\sum_{j=1}^n\sigma_{ij}(t)\,dW_j(t)\right],\quad0\leq
t\leq\infty,\quad i=1,\dots,d.
$$
Here $W=(W_{1},\cdots,W_{n})^{T}$ is a $n-$dimensional Brownian motion on a filtered probability space
$(\Omega,\{\sF_t\}_{0\leq t\leq T},\mathcal{F},\mathbb{P}),$ where $\{\sF_t\}_{0\leq t\leq T}$ is the completed filtration generated by $W.$ Here we assume $d\leq n$, i.e., there are at least as many sources
of uncertainty as assets. As usual $\{\alpha(t)\}_{t\in[0,\infty)}=\{(\alpha_i(t))_{i=1,\cdots,d}\}_{t\in[0,\infty)}$ is an $\mathbb{R}^{d}$ valued \textit{mean rate of return} process,
and\newline $\{\sigma(t)\}_{t\in[0,\infty)}=\{(\sigma_{ij}(t))_{i=1,\cdots,d}^{j=1,\cdots,n}\}_{t\in[0,\infty)}$ is an $d\times n-$matrix valued \textit{volatility} process,
 progressively measurable with respect to $\{\sF_t\}_{0\leq t\leq T}.$  

\textit{Standing Assumption 2.1}
The process $\alpha$ is uniformly bounded and the volatility matrix $\sigma$ has full rank. Moreover $\sigma\sigma^{T}$ is assumed uniformly elliptic, i.e. $K I_{d}\geq \sigma\sigma^{T}\geq \epsilon I_{d},$ for some $K\geq\epsilon>0.$ 

This implies the existence of the inverse $(\sigma(t)\sigma^\mathrm{T}(t))^{-1}$ and the market price of risk process
\begin{equation}\label{1}\tilde{\theta}(t)=\sigma^\mathrm{T}(t)(\sigma(t)\sigma^\mathrm{T}(t))^{-1}\alpha(t), \end{equation}
which is uniformly bounded. All the processes encountered are defined on the fixed, finite interval $[0,T].$ The stochastic exponential process 
 
  \begin{equation}\label{Nov}
  \tilde{Z}(t)=Z_{\tilde{\theta}}(t)\triangleq\exp\left\{-\int_{0}^{t}\tilde{\theta}^{T}(u)\,dW(u)-\frac{1}{2}\int_{0}^{t}\parallel\tilde{\theta}(u)\parallel^{2}\,du\right\}\end{equation}
is a (true) martingale, thus by the Girsanov theorem (section $3.5$ in \cite{KAR}) \begin{equation}\label{Br}\tilde{W}(t)=W(t)+\int_{0}^{t}\tilde{\theta}(u)\,du \end{equation} is a Brownian motion under the equivalent martingale measure
  \begin{equation}\label{Q}
  \mathbb{\tilde{Q}}(A)\triangleq \mathbb{E}[\tilde{Z}(T){\bf{1}}_{A}],\qquad A\in \sF_T.\end{equation}
  
  \begin{definition}
  We denote by $\mathcal{M}$ the set of probability measure satisfying
\begin{description}
\item[(i)]
$\Q\ll\P$ and $\frac{d\Q}{d\P}\in L_T^2(\P)$;
\item[(ii)]
$S$ is a local martingale under $\Q$
on $[0,T].$ 
\end{description}
  \end{definition}
 In the light of boundedness of $\tilde{\theta},$ it follows that  $\mathbb{\tilde{Q}}\in\mathcal{M}.$   
   
\subsection{Portfolio and wealth processes}

 A (self-financing) portfolio is defined as a pair $(x,\pi).$ The constant $x,$ exogenously given, is the initial value of the
 portfolio and $\pi=(\pi_{1},\cdots,\pi{d})^{T}$ is a predictable $S-$ integrable process which specifies how many units of the asset
 $i$ are held in the portfolio at time $t.$ The wealth process of such a portfolio is given by
 \begin{equation}\label{welath}
 X^{x,\pi}(t)=x+\int_{0}^{t}\pi(u)^{T}dS(u).
 \end{equation}

\subsection{Utility Function}
 A function $U:(0,\infty)\rightarrow \mathbb{R}$ strictly increasing and strictly concave is called utility function.
 We restrict  ourselves to utility functions which are $3-$times continuous differentiable and satisfy the Inada conditions
 \begin{equation}\label{In}
 U'(0+)\triangleq\lim_{x\downarrow 0}U'(x)=\infty,\quad U'(\infty)\triangleq\lim_{x\uparrow{\infty}}U'(x)=0.
 \end{equation}
We shall denote by $I(\cdot)$ the (continuous, strictly decreasing) inverse of the marginal utility function $U'(\cdot),$ and by \eqref{In}
\begin{equation}\label{In1}
 I(0+)\triangleq\lim_{x\downarrow 0}I(x)=\infty,\quad I(\infty)\triangleq\lim_{x\uparrow{\infty}}I(x)=0.
 \end{equation}
Let us introduce the Legeandre transform of $-U(-x)$
\begin{equation}\label{Leg}
\tilde{U}(y)\triangleq\sup_{x>0}[U(x)-xy]=U(I(y))-yI(y),\quad 0<y<\infty.
\end{equation}
The function $\tilde{U}(\cdot)$ is strictly decreasing, strictly convex and satisfies  the dual relationships
\begin{equation}\label{In12}
\tilde{U}'(y)=-x\quad\mbox{iff}\quad U'(x)=-y,
 \end{equation}
 and
 \begin{equation}\label{In13}
U(x)=\inf_{y>0}[\tilde{U}(y)+xy]=\tilde{U}(U'(x))+xU'(x),\quad 0<x<\infty.
 \end{equation}
\textit{Standing Assumption $2.1$}
\begin{equation}\label{0op}
y^{2}|I''(y)|\vee(-yI'(y))\vee I(y)<k_1y^{-\kappa}\quad\mbox{for\,\,every}\,\,y\in(0,\infty),
\end{equation}
for some $k_1>0,$ $\kappa>0,$ and $a\vee b=\max(a,b).$

\section{Objective}
Let $x$ be the agent's initial wealth and $V_{x}$ a process progressively measurable with respect to $\{\sF_t\}_{0\leq t\leq T}$ which satisfies $V_{x}(0)=0,$ 
\begin{equation}\label{as100}
\mathbb{E}[\sup_{0\leq t\leq T}|V_{x}(t)|^{2}]<\infty.
\end{equation}  
For a given initial positive wealth $x$ and a given utility function $U$ find a process $V_{x}$ as above
such that for any stopping time $\tau\leq T$ there exists a portfolio process $\hat{\pi}_{x}$ such that
\begin{equation}\label{eq3}
{\sup_{\pi\in\mathcal{A}_{V}(x,\tau)}}\mathbb{E}U(X^{x,\pi}(\tau)+V_{x}(\tau))=\mathbb{E}U(X^{x,\hat{\pi}}(\tau)+V_{x}(\tau)).
\end{equation} 
Here $\mathcal{A}_{V}(x,\tau)$ is the set of admissible portfolios given the time horizon $\tau.$ It is defined by

\begin{eqnarray*} \mathcal{A}_{V}(x,\tau)&\triangleq&\left\{\pi \left|
\begin{array}{l}
X^{x,\pi}(\tau)+V_{x}(\tau)>0,\,\,\mathbb{E}[U(X^{x,\pi}(\tau)+V_{x}(\tau))]^{-}<\infty,\,\, \mbox{and}\\
 \,X^{x,\pi}(t)\,\,\,\mbox{is\,\,$\mathbb{Q}$-\,supermartingale,}\,\,0\leq t\leq\tau,\,\forall\,\,\mathbb{Q}\in\mathcal{M}
\end{array}\right.\right\}.
\end{eqnarray*}

\begin{remark}
Notice this class of admissible trading strategy includes the classical one
\begin{equation}\label{leq4}
\mathcal{A}^{*}_{V}(x,\tau)\triangleq\{\pi|\,\,X^{x,\pi}(t)+V_{x}(t)>0,\,\,0\leq t\leq\tau\}.
\end{equation}

 \noindent Indeed the process $X^{x,\pi}(t)$ is a local martingale under the probability measure $\mathbb{Q}.$
 Let $T_{n}\uparrow T$ be a localizing sequence, $u\leq t\leq\tau,$ and $\mathbb{E}^{\mathbb{Q}}$ the expectation operator
 with respect to a given martingale measure $\mathbb{Q}$. H\"older's inequality and \eqref{as100} imply $\mathbb{E}^{\mathbb{Q}}[\sup_{0\leq t\leq T}|V_{x}(t)|]<\infty.$ Fatou's Lemma and \eqref{as100} yield
 
 \begin{eqnarray*}{\mathbb{E}^{\mathbb{Q}}}[X^{x,\pi}(t)+\sup_{0\leq s\leq t}|V_{x}(s)|\,\,|\sF_u ]&\leq &{\liminf_{n\rightarrow\infty}}\,{\mathbb{E}^{\mathbb{Q}}}[X^{x,\pi}(t\wedge T_{n})+\sup_{0\leq s\leq t}|V_{x}(s)|\,\,|\sF_u]\\&=&X^{x,\pi}(u)+{\mathbb{E}^{\mathbb{Q}}}[\sup_{0\leq s\leq t}|V_{x}(s)|\,\,|\sF_u] , \end{eqnarray*}
whence 
\begin{equation}\label{qwe009}
{\mathbb{E}^{\mathbb{Q}}}[X^{x,\pi}(t)\,\,|\sF_u ]\leq X^{x,\pi}(u),
\end{equation}   
so 
$X^{x,\pi}(t)$\,\,\,is\,\,\,$\mathbb{Q}$-\,supermartingale,\,\,$0\leq t\leq\tau.$

\begin{flushright}
$\diamond$
\end{flushright}

\end{remark}

\section{The complete market solution}

In this section we prove the existence and uniqueness of the process $V_{x}$ if $d=n$ (complete market).
It will be helpful first to solve the utility maximization of the corrected final wealth at some stopping
time horizon $\tau,$ i.e.,

\begin{equation}\label{p00}
{\sup_{\xi_{\tau}\in\mathcal{B}_{V}(\tau,x)}}\mathbb{E}U(\xi_{\tau}),
\end{equation}
where
\begin{equation}\label{eq40}
\mathcal{B}_{V}(\tau,x)\triangleq\{\xi_{\tau}:\sF_\tau\,\,\mbox{measurable,\,positive},\,\,\tilde{\mathbb{E}}\xi_{\tau}\leq x+\tilde{\mathbb{E}}V_{x}(\tau),\,\,\,\mathbb{E}[U(\xi_{\tau})]^{-}<\infty\}.
\end{equation}

 \noindent From now on we fix the stopping time $\tau.$ According to \eqref{0op} 
$$\mathbb{E}\,{\sup_{0\leq t\leq T}}\tilde{Z}(t)I(\lambda\tilde{Z}(t))\leq k_{1} \mathbb{E}\,{\sup_{0\leq t\leq T}}\tilde{Z}(t)^{1-\kappa}.  $$ If $\kappa\leq 1$ Burholder-Davis-Gundy (see p. $166$ in \cite{KAR}) in conjunction with boundedness of $\tilde{\theta}$ and  H\"older's inequality yields 
 $$\mathbb{E}\,{\sup_{0\leq t\leq T}}\tilde{Z}(t)^{1-\kappa}<\infty. $$ Otherwise $$\tilde{Z}(t)^{1-\kappa}=M(t)^{\kappa-1}\exp\left\{\int_{0}^{t}(\kappa-1)\parallel\tilde{\theta}(u)\parallel^{2}\,du\right\},$$
where $M(t)$ is the local martingale
$$M(t)\triangleq \exp\left\{\int_{0}^{t}\tilde{\theta}^{T}(u)\,dW(u)-\frac{1}{2}\int_{0}^{t}\parallel\tilde{\theta}(u)\parallel^{2}\,du\right\}, $$ and the same argument applies. Therefore 
\begin{equation}\label{kl}\mathbb{E}\,{\sup_{0\leq t\leq T}}\tilde{Z}(t)I(\lambda\tilde{Z}(t))<\infty.\end{equation}

\noindent Similarly by \eqref{as100}
$$\mathbb{E}\,{\sup_{0\leq t\leq T}}\tilde{Z}(t)|V_{x}(t)|<\infty, $$
 whence
 \begin{equation}\label{01p}
\mathbb{E}\,{\sup_{0\leq t\leq T}}\tilde{Z}(t)[I(\lambda\tilde{Z}(t))-V_{x}(t)]<\infty,
\end{equation}
for every $\lambda>0.$ Therefore the function ${\mathcal{X}_{\tau}}:(0,\infty)\rightarrow (-\mathbb{E}[\tilde{Z}(\tau)V_{x}(\tau)],\infty)$ given by
\begin{equation}\label{11p}
\mathcal{X}_{\tau}(\lambda)\triangleq\mathbb{E}\tilde{Z}(\tau)[I(\lambda\tilde{Z}(\tau))-V_{x}(\tau)],
\end{equation} 
 is continuous and strictly decreasing. We shall denote by $\mathcal{Y}_{\tau}(\cdot)$ its continuous, strictly decreasing inverse.

\begin{lemma}\label{LL:3}
The random variable 
\begin{equation}\label{9i}
\hat{\xi}_{\tau}\triangleq I(\mathcal{Y}_{\tau}(x)\tilde{Z}(\tau)),
\end{equation} 
 belongs to $\mathcal{B}_{V}(\tau,x).$
Moreover for any $\xi_{\tau}\in\mathcal{B}_{V}(\tau,x)$
\begin{equation}\label{29i}
\mathbb{E}U(\xi_{\tau})\leq \mathbb{E}U(\hat{\xi}_{\tau}),
\end{equation}
and $\hat{\xi}_{\tau}$ is unique with this property. 
\end{lemma}
\noindent {\sc Proof:} See the Appendix. 
\begin{flushright}
$\square$
\end{flushright}  

\noindent The above Lemma solves the optimization problem at the level of claims; to obtain the optimal portfolio we proceed as follows.
 The martingale $\tilde{\mathbb{E}}[(\hat{\xi}_{\tau}-V_{x}(\tau))|\sF_t]$ admits the stochastic integral representation
 \begin{equation}\label{18i}
\tilde{\mathbb{E}}[(\hat{\xi}_{\tau}-V_{x}(\tau))|\sF_t]=x+\int_{0}^{t}\psi^{T}(u)d\tilde{W}(u),\qquad 0\leq t\leq\tau,
\end{equation}
 for some $\sF_t-$adapted process $\psi(\cdot)$ that satisfies $\int_{0}^{\tau} ||\psi(u)||^{2}\,du<\infty$ a.s. (e.g., \cite{KAR1}, Lemma $1.6.7$).
 Let
 $$X^{x,\hat{\pi}}(t)\triangleq\tilde{\mathbb{E}}[(\hat{\xi}_{\tau}-V_{x}(\tau))|\sF_t],$$
 \begin{equation}\label{4r}
 \hat{\pi}_{x}(t)=(\sigma^{T}(t))^{-1}\psi(t) \qquad 0\leq t\leq\tau,
 \end{equation}
 
\noindent and notice that
 \begin{equation}\label{eu}
X^{x,\hat{\pi}}(\tau)+V_{x}(\tau)=\hat{\xi}_{\tau}\triangleq I(\mathcal{Y}_{\tau}(x)\tilde{Z}(\tau)).
 \end{equation}
For any $\pi\in\mathcal{A}_{V}(x,\tau)$ the corrected final wealth $X^{x,\pi}(\tau)+V_{x}(\tau)\in\mathcal{B}_{V}(\tau,x)$ and by Lemma \ref{LL:3}, 
 \begin{equation}\label{0eq3}
{\sup_{\pi\in\mathcal{A}_{V}(x,\tau)}}\mathbb{E}U(X^{x,\pi}(\tau)+V_{x}(\tau))=\mathbb{E}U(X^{x,\hat{\pi}}(\tau)+V_{x}(\tau)).
\end{equation}

\begin{theorem}\label{main1}
There exists a unique process $V_{x}$ which satisfies $V_{x}(0)=0$ and \eqref{as100}, for which there exists a portfolio process $\hat{\pi}_{x}(\cdot)$ such that for any stopping time $\tau,$
\begin{equation}\label{eq3}
{\sup_{\pi\in\mathcal{A}_{V}(x,\tau)}}\mathbb{E}U(X^{x,\pi}(\tau)+V_{x}(\tau))=\mathbb{E}U(X^{x,\hat{\pi}}(\tau)+V_{x}(\tau)).
\end{equation} 
It is given by
\begin{equation}\label{q3}
V_{x}(t)=\int_{0}^{t}F(U'(x)\tilde{Z}(u))||\tilde{\theta}(u)||^{2}\,du,\quad 0\leq t\leq T,
\end{equation}
where
\begin{equation}\label{0q3}
F(z)=\frac{1}{2}I''(z)z^{2}+I'(z)z.
\end{equation}
The optimal portfolio process $\hat{\pi}_{x}(\cdot)$ is
\begin{equation}\label{q34}
(\hat{\pi}_{x}(t))_{i}=-\frac{1}{S_{i}(t)}((\sigma^{T}(t))^{-1}I'(U'(x)\tilde{Z}(t))U'(x)\tilde{Z}(t)\tilde{\theta}(t))_{i}, \quad 0\leq t\leq T,\,\, i=1,\cdots,d. 
\end{equation}

\end{theorem} 

\noindent The proof relies on the following Lemma.

\vspace{0.3cm}

\begin{lemma}\label{L:44}
Let $\tilde{V}_{x}$ be a process satisfying $\tilde{V}_{x}(0)=0$ and \eqref{as100}. There exists a portfolio process $\tilde{\pi}(\cdot)$ such that for any stopping time $\tau,$
\begin{equation}\label{eq3}
{\sup_{\pi\in\mathcal{A}_{\tilde{V}}(x,\tau)}}\mathbb{E}U(X^{x,\pi}(\tau)+\tilde{V}_{x}(\tau))=\mathbb{E}U(X^{x,\tilde{\pi}}(\tau)+\tilde{V}_{x}(\tau)),
\end{equation} 
if and only if 
\begin{equation}\label{1eq3}
X^{x,\tilde{\pi}}(t)+\tilde{V}_{x}(t)=I(U'(x)\tilde{Z}(t)),\quad 0\leq t\leq T.
\end{equation}

\end{lemma}

\noindent {\sc Proof:} See the Appendix.

\begin{flushright}
$\square$
\end{flushright} 

\noindent {\sc Proof of Theorem \ref{main1}:} It\^{o}'s Lemma gives
\begin{eqnarray}\label{kg}
I(U'(x)\tilde{Z}(t))&=&x-\int_{0}^{t}U'(x)I'(U'(x)\tilde{Z}(u))\tilde{Z}(u)\tilde{\theta}(u)d\tilde{W}(u)\\\notag & & \mbox{}+\int_{0}^{t}F(U'(x)\tilde{Z}(u))||\tilde{\theta}(u)||^{2}\,du.
\end{eqnarray}
For any $\pi,$ this and \eqref{welath} yield
\begin{eqnarray}\notag
\tilde{V}_{x}(t)&\triangleq&I(U'(x)\tilde{Z}(t))-X^{x,\pi}(t)\\\label{lo}&=&\int_{0}^{t}\bar{\pi}(u)dS(u)+\int_{0}^{t}F(U'(x)\tilde{Z}(u))||\tilde{\theta}(u)||^{2}\,du,
\end{eqnarray}
for the process $$\bar{\pi}_{i}(t)=-\left[\pi_{i}(t)+\frac{1}{S_{i}(t)}((\sigma^{T}(t))^{-1}I'(U'(x)\tilde{Z}(t))\tilde{Z}(t)\tilde{\theta}(t))_{i}\right]=-[\pi_{i}(t)-\hat{\pi}_{i}(t)],$$
with $\hat{\pi}_{x}$ of \eqref{q34}. If we set $\pi=\hat{\pi}_{x},$ then $\bar{\pi}$ is the zero vector so $\tilde{V}_{x}(t)=V_{x}(t)$ of \eqref{q3}, and \eqref{1eq3} holds. Therefore by Lemma \ref{L:44}
\begin{equation}\label{eq33}
{\sup_{\pi\in\mathcal{A}_{V}(x,\tau)}}\mathbb{E}U(X^{x,\pi}(\tau)+V_{x}(\tau))=\mathbb{E}U(X^{x,\hat{\pi}}(\tau)+V_{x}(\tau)),
\end{equation}
with $V_{x}$ and $\hat{\pi}_{x}$ defined by \eqref{q3} and \eqref{q34}. From \eqref{kg} and \eqref{lo} follows
\begin{equation}\label{win}
X^{x,\hat{\pi}}(t)=x-\int_{0}^{t}U'(x)I(U'(x)\tilde{Z}(u))\tilde{Z}(u)\tilde{\theta}(u)d\tilde{W}(u).
\end{equation}
 The inequality
$$U(I(y))\geq U(x)+y[I(y)-x]\qquad \forall x>0,\,y>0,$$
follows from \eqref{Leg} and in conjunction with \eqref{0op} give
\begin{equation}\label{tris}
U^{-}(I(y))<k_{2}+k_{3}y^{1-\kappa},
\end{equation}
for some positive $k_{2},$ $k_{3}.$ This and the boundedness of $\tilde{\theta}$ imply that
for any $\tau$ stopping time 
$$\mathbb{E}[U(X^{x,\hat{\pi}}(\tau)+V_{x}(\tau))]^{-}=\mathbb{E}U[I(U'(x)\tilde{Z}(\tau))]^{-}<\infty.$$ Moreover by Remark $3.1,$ $\hat{\pi}_{x}\in\mathcal{A}_{V}(x,\tau).$
 The assumption \eqref{0op} show that the process $V_{x}$ of \eqref{q3} satisfies \eqref{as100}. The uniqueness of $V_{x}$ is up to translations by $S$-integrals, i.e., the process $\tilde{V}_{x}(t)=\hat{V}_{x}(t)+\int_{0}^{t}\bar{\pi}(u)dS(u)$ is another solution and the corresponding optimal portfolio is $\tilde{\pi}=\hat{\pi}_{x}-\bar{\pi}.$

\begin{flushright}
$\square$
\end{flushright}

\section{The incomplete market solution}
When markets are incomplete we establish existence of the process $V_{x},$ but we can no longer prove uniqueness.
\begin{theorem}\label{main2}
There exists a process $V_{x}$ which satisfies $V_{x}(0)=0,$ and \eqref{as100}, for which there exists a portfolio process $\hat{\pi}_{x}(\cdot)$ such that for any stopping time $\tau,$
\begin{equation}\label{eq4}
{\sup_{\pi\in\mathcal{A}_{V}(x,\tau)}}\mathbb{E}U(X^{x,\pi}(\tau)+V_{x}(\tau))=\mathbb{E}U(X^{x,\hat{\pi}}(\tau)+V_{x}(\tau)).
\end{equation} 
It is given by
\begin{equation}\label{q4}
V_{x}(t)=\int_{0}^{t}F(U'(x)\tilde{Z}(u))||\tilde{\theta}(u)||^{2}\,du,\quad 0\leq t\leq T,
\end{equation}
where $F$ is the function of \eqref{0q3}. The optimal portfolio $\hat{\pi}_{x}$ is
\begin{equation}\label{q35}
(\hat{\pi}_{x}(t))_{i}=-\frac{1}{S_{i}(t)}((\sigma(t)\sigma^{T}(t))^{-1}I'(U'(x)\tilde{Z}(t))\tilde{Z}(t)\alpha(t))_{i}, \quad 0\leq t\leq T,\,\, i=1,\cdots,d, 
\end{equation}
where $\alpha$ is the \textit{mean rate of return} process.
\end{theorem}
\noindent {\sc Proof:} It follows from  It\^{o}'s Lemma applied to $I(U'(x)\tilde{Z}(t))$ that 
\begin{equation}\label{1eq4}
X^{x,\hat{\pi}}(t)+V_{x}(t)=I(U'(x)\tilde{Z}(t)),\quad 0\leq t\leq T,
\end{equation}
with $V_{x}$ defined in \eqref{q4} and $\hat{\pi}_{x}$ of \eqref{q35}. The same argument as in the proof of Lemma \ref{L:44} concludes.

\begin{flushright}
$\square$
\end{flushright}

\begin{remark}
If the utility is logarithmic, i.e.,
$U(x)=\log{x},$ then $V_{x}(t)\equiv 0$. Being optimal for the logarithmic utility, the vector $\hat{\pi}_{x}$ satisfies
\begin{equation}\label{0i}(\hat{\pi}_{x}(t))_{i}=\frac{(\zeta_{M}(t))_{i}X^{\hat{\pi}}(t)}{S_{i}(t)},\quad i=1,\dots,d,\end{equation} with $\zeta_{M}(t)\triangleq(\sigma(t)\sigma^\mathrm{T}(t))^{-1}\mu(t)$ the Merton proportion. The future evolution of $S$ does not enter in the formula (\ref{q35}) and (\ref{0i}), hence we refer to $\hat{\pi}_{x}$  as \textit{the myopic component}.
\end{remark}

\begin{remark}
Direct computation shows that $F$ of \eqref{0q3} is in fact
$$F(y)=\frac{y^{2}}{[U''(I(y))]^{2}}\left[-\frac{U'''(I(y))}{U''(I(y))}+2\frac{U''(I(y))}{U'(I(y))} \right], $$
and it has the following economic interpretation. It is the difference between the coefficient of prudence,$-\frac{U'''(I(y))}{U''(I(y))},$ and twice the coefficient
of relative risk aversion, $-\frac{U''(I(y))}{U'(I(y))}.$ The coefficient of prudence reflects an individual's propensity to take precautions
when faced with risk. The coefficient of relative risk aversion goes back to Arrow-Pratt and reflects the tendency to avoid risk altogether.
\end{remark}

\section{The hedging portfolio}
We have learnt from the previous sections that we have to top up investor's wealth with the process $V_{x}$ of \eqref{q4} in order to achieve time consistency
of the optimal investment. Let us compute the process $V_{x}$ for different utility functions using \eqref{q4}. 
In the case of an exponential utility, i.e., $U(x)=-e^{-ax}$:
$$V_{x}(t)=\frac{1}{a}\int_{0}^{t}||\tilde{\theta}(u)||^{2}\,du. $$
If the utility is CRRA, i.e., $U(x)=\frac{x^{p}}{p}$:
$$V_{x}(t)=\frac{xp}{2(p-1)^{2}}\int_{0}^{t}[\tilde{Z}(u)]^{\frac{1}{p-1}}||\tilde{\theta}(u)||^{2}\,du.$$
This shows that for different utility functions the process $V_{x}$ can be either positive or negative. In the general case it can be decomposed as the difference of two 
increasing processes $V_{x}(t)=V_{x}^{+}(t)-V^{-}_{x}(t),$ as follows:
\begin{equation}\label{ln}
V^{+}_{x}(t)= \int_{0}^{t}F^{+}(U'(x)\tilde{Z}(u))||\tilde{\theta}(u)||^{2}\,du,\quad 0\leq t\leq T,
\end{equation}

\begin{equation}\label{ln1}
V^{-}_{x}(t)= \int_{0}^{t}F^{-}(U'(x)\tilde{Z}(u))||\tilde{\theta}(u)||^{2}\,du,\quad 0\leq t\leq T,
\end{equation}
where as usual $a^{-}\triangleq\max\{-a,0\},$ and $a^{+}\triangleq\max\{a,0\}.$
As we have already seen the agent has to have her wealth adjusted by the process $V_{x}$ at all times. The natural question is how can she
achieve this ? We answer this question under the assumption that she starts with some initial wealth, she does not receive extra funds in the future
 and the only investment instruments are the stocks and the bond. In order to hold $-\{V^{-}_{x}(T)\}_{0\leq t\leq T}$ in her portfolio the agent should
 consume at rate $\{F^{-}(U'(x)\tilde{Z}(t))||\tilde{\theta}(t)||^{2}\}_{0\leq t\leq T}$. The process $\{V^{+}_{x}(T)\}_{0\leq t\leq T}$ is positive and
 increasing and in order to finance it, the agent should allocate some of her initial wealth to its generation. One approach is to create a portfolio, which at the final time
  $T$ replicates $V^{+}_{x}(T).$ The value of such portfolio at any time $t$ will exceed $V_{x}^{+}(t)$ and the agent can consume the difference. The drawback of
  this methodology is that incomplete markets makes perfect hedging impossible, hence the replication of $V_{x}^{+}(T)$ is impossible. The resolution we propose
  is to consider a \textit{risk-minimizing} strategy as in \cite{FO}. Strictly speaking it is the strategy for which the remaining risk (in hedging the contingent claim $V_{x}^{+}(T)$ ) is minimal under all infinitesimal perturbations of the strategy at some intertemporal time t (see \cite{FO}). It can be determined using the Kunita-Watanabe projection technique and
  is intimately related to \textit{minimal martingale measure} which in our setup is $\mathbb{\tilde{Q}}$. Let us consider the process $\{\bar{V}_{x}^{+}(t)\}_{0\leq t\leq T},$ $\bar{V}_{x}^{+}(t)\triangleq\tilde{\mathbb{E}}[V_{x}^{+}(T)|\sF_{t}].$ By Corollary $1$ page $181$ in \cite{Protter} or Proposition $4.14$ page $181$ in \cite{KAR} one gets the Kunita-Watanabe decomposition of $\{\bar{V}_{x}^{+}(t)\}_{0\leq t\leq T},$   
\begin{equation}\label{jh}\bar{V}_{x}^{+}(t)=\int_{0}^{t}\bar{\pi}_{x}(u)dS(u)+M(t)+\tilde{\mathbb{E}}V_{x}^{+}(T),\end{equation}
 where $M$ is a process orthogonal to $S,$ i.e., $\langle M, S\rangle=0$ and $M(0)=0.$ Moreover the process $M$ is a martingale under both $\mathbb{P}$ and $\mathbb{Q}.$ According to \cite{FO}, $\bar{\pi}_{x}$ is the \textit{risk-minimizing} strategy and we refer to it as \textit{the hedging portfolio}. The initial cost to implement it is $\tilde{\mathbb{E}}V_{x}^{+}(T),$ and requires some of the initial wealth.
 Let 
$$x_{*}=\inf\{z\geq 0|z+\tilde{\mathbb{E}}V_{z}^{+}(T)=x\},$$ which by \eqref{0op} exists and is positive. 

\begin{theorem}\label{main3}
Starting with the initial wealth $x,$ the agent should invest $\hat{\pi}_{x_{*}}+\bar{\pi}_{x_{*}}$ in stocks, where
$$(\hat{\pi}_{x_{*}}(t))_{i}=-\frac{1}{S_{i}(t)}((\sigma(t)\sigma^{T}(t))^{-1}I'(U'(x_{*})\tilde{Z}(t))\tilde{Z}(t)\alpha(t))_{i}, \quad 0\leq t\leq T,\,\, i=1,\cdots,d,$$
is \textit{the myopic portfolio}, and $\bar{\pi}_{x_{*}}$ of \eqref{jh} is \textit{the hedging portfolio}. She can consume $V^{-}_{x_{*}}(t)+\tilde{\mathbb{E}}[V_{x_{*}}^{+}(T)-V_{x_{*}}^{+}(t)|\sF_{t}].$ This investment strategy is time consistent up to the \textit{intrinsic} risk $M$ of \eqref{jh}.  
\end{theorem}
\noindent {\sc Proof:} At any time $t,$ up to $M(t)$ the agent's wealth following the above policy is
$$X^{x_{*},\hat{\pi}}(t)+V_{x_{*}}(t),$$
so Theorem \ref{main2} concludes. 
\begin{flushright}
$\square$
\end{flushright}

\noindent In the reminder of this section we show how to compute $\bar{\pi}_{x_{*}}$ explicitly. The martingale $\bar{V}^{+}_{x_{*}}(t)\triangleq\tilde{\mathbb{E}}[V_{x_{*}}^{+}(T)|\sF_{t}],$\,\,admits the stochastic integral representation
 \begin{equation}\label{8i}
\bar{V}^{+}_{x_{*}}(t)=\tilde{\mathbb{E}}V_{x_{*}}^{+}(T)+\int_{0}^{t}\beta^{T}(u)d\tilde{W}(u),\qquad 0\leq t\leq T,
\end{equation}
 for some $\sF_t-$adapted process $\beta(\cdot)$ which satisfies $\int_{0}^{T} ||\beta(u)||^{2}\,du<\infty$ a.s. (e.g., \cite{KAR1}, Lemma $1.6.7$). In light of this we want $$\int_{0}^{t}\beta(u)d\tilde{W}(u)=\int_{0}^{t}\bar{\pi}_{x_{*}}(u)dS(u)+M(t),\quad 0\leq t\leq T.$$ We are looking for the orthogonal process $M$
of the form $M=\delta\cdot\tilde{W}$. Let $A$ be the $d\times n$ matrix with the entries $A_{ij}=S_i\sigma_{ij}.$ Since the volatility matrix $\sigma$ has linearly independent rows, the matrix $A$ has linearly independent rows, i.e., $\mbox{rank}\,A=d.$
The processes $\bar{\pi}_{x_{*}}$ and $\delta$ must solve
\begin{equation}\label{g}
A^{T}\bar{\pi}_{x_{*}}+\delta=\beta,\quad\mbox{and}\quad A\delta=0,
\end{equation}
where the second equality comes from 
$$ \langle M,S_i \rangle=\int_{0}^{\cdot}\sum_{k=1}^{n}\delta_{k}A_{ik}\,dt=\int_{0}^{\cdot}(A\delta)_{i}\,dt =0.   $$

\noindent Since $\mathrm{Im}(A^{T})\oplus\mathrm{Ker}(A)=\mathbb{R}^{n},$ we can uniquely decomposed $\beta=z_1+z_2,$ with $z_1\in \mbox{Im}(A^{T})$ and $z_2\in\mathrm{Ker}(A).$ Set $\delta=z_{2}$ and $\bar{\nu}_{x}=y,$ where $z_1=A^{T}y$ ( the uniqueness
of $\bar{\pi}_{x_{*}}$ is due to $\mbox{rank}\,A=d$). Hence
\begin{equation}\label{o}
\bar{V}^{+}_{x_{*}}=\bar{\pi}_{x_{*}}\cdot S+M+\tilde{\mathbb{E}}V_{x_{*}}^{+}(T)
,\end{equation}with

\begin{equation}\label{x}
\bar{\pi}_{x_{*}}=y,\quad\mbox{and}\quad M=\int z_{2}\,d\tilde{W}.
\end{equation}Let us notice that the representation formula in (\ref{o})
takes place under the probability measure $\tilde{\mathbb{Q}}.$ However by Theorem $14$ page $60$ in \cite{Protter} this takes place under $\mathbb{P},$ since $\mathbb{P}\sim \tilde{\mathbb{Q}}.$ Moreover $\langle M, S\rangle=0$ and $M(0)=0$ under $\mathbb{P}.$ The process $\beta(t)$ of \eqref{8i} can be computed explicitly by Clark-Haussmann-Ocone formula. 

\section{Appendix}

\noindent {\bf {Proof of Lemma \ref{LL:3}}}:
Let us notice that
\begin{equation}\label{i8i}
\mathbb{E}[\tilde{Z}(\tau)(\hat{\xi}_{\tau}-V_{x}(\tau))]=\mathcal{X}_{\tau}(\mathcal{Y}_{\tau}(x))=x.
\end{equation}
  
 \noindent Boundedness of $\tilde{\theta}$ and \eqref{tris} give
 $$\mathbb{E}[U(\hat{\xi}_{\tau})]^{-}<\infty,$$
  whence $\hat{\xi}_{\tau}\in\mathcal{B}_{V}(\tau,x).$ Let $\xi_{\tau}\in\mathcal{B}_{V}(\tau,x),$ then  by \eqref{i8i}
 $$\mathbb{E}\tilde{Z}(\tau)\xi_{\tau}\leq x+\tilde{\mathbb{E}}V_{x}(\tau)= \mathbb{E}\tilde{Z}(\tau)\hat{\xi}_{\tau}.$$ This and the inequality \eqref{In13} give
 \begin{eqnarray*}\mathbb{E}U(\xi_{\tau})&\leq&\mathbb{E}\tilde{U}(\mathcal{Y}_{\tau}(x)\tilde{Z}(\tau))+\mathcal{Y}_{\tau}(x)\,\mathbb{E}\tilde{Z}(\tau)\xi_{\tau}\\
 &\leq&\mathbb{E}\tilde{U}(\mathcal{Y}_{\tau}(x)\tilde{Z}(\tau))+\mathcal{Y}_{\tau}(x)\,\mathbb{E}\tilde{Z}(\tau)\hat{\xi}_{\tau} =\mathbb{E}U(\hat{\xi}_{\tau}), \end{eqnarray*}
for any $\xi_{\tau}\in\mathcal{B}_{V}(\tau,x).$ Uniqueness of $\hat{\xi}_{\tau}$ is a consequence of the concavity of $U.$
 
\begin{flushright}
$\diamond$
\end{flushright}

\noindent {\bf{Proof of Lemma \ref{L:44}}}: Let us first establish the sufficiency. For any stopping time $\tau,$ and $\pi\in\mathcal{A}_{V}(x,\tau),$
the process $X^{x,\pi}$ is a $\mathbb{\tilde{Q}}$ supermartingale and Problem $3.26,$ p. $20$ in \cite{KAR} yields  
\begin{equation}\label{i9}
\mathbb{E}[\tilde{Z}(\tau)X^{x,\pi}(\tau)]\leq x.
\end{equation} 
Recall that
$$X^{x,\tilde{\pi}}(t)=I(U'(x)\tilde{Z}(t))-\tilde{V}_{x}(t),$$
cf \eqref{1eq3}. The process $X^{x,\tilde{\pi}},$ being a stochastic integral with integrator $S,$
is a local $\mathbb{Q}-$martingale. Arguing as in \eqref{kl} one gets $\mathbb{E}\,{\sup_{0\leq t\leq T}}|I(U'(x)\tilde{Z}(t))|^{2}<\infty,$ and hence by \eqref{as100}\newline $\mathbb{E}\,{\sup_{0\leq t\leq T}}|X^{x,\tilde{\pi}}(t)|^{2}<\infty.$ H\"older's inequality implies $\tilde{\mathbb{E}}\,{\sup_{0\leq t\leq T}}|X^{x,\tilde{\pi}}(t)|<\infty.$ Since for any $0\leq t\leq T,$ $X^{x,\tilde{\pi}}(t)\leq {\sup_{0\leq t\leq T}}|X^{x,\tilde{\pi}}(t)|,$ by Dominated Convergence Theorem, the process $X^{x,\tilde{\pi}}$ is a (true) $\mathbb{Q}-$martingale.  
 Hence for any stopping time $\tau$ 
\begin{equation}\label{1i9}
\mathbb{E}[\tilde{Z}(\tau)X^{x,\tilde{\pi}}(\tau)]=x,
\end{equation}  
by Problem $3.26,$ p. $20$ in \cite{KAR}. Due to concavity of utility function

  \begin{eqnarray*}
  U(X^{x,\pi}(\tau)+V_{x}(\tau))-U(X^{x,\tilde{\pi}}(\tau)+V_{x}(\tau))&\leq& U'(X^{x,\tilde{\pi}}(\tau)+V_{x}(\tau))(X^{x,\pi}(\tau)-X^{x,\tilde{\pi}}(\tau))\\&=&U'(x)\tilde{Z}(\tau)(X^{x,\pi}(\tau)-X^{x,\tilde{\pi}}(\tau)).  \end{eqnarray*}
Taking expectation and using \eqref{i9} and \eqref{1i9} we get the sufficiency.

As for necessity of \eqref{1eq3}, let us assume the existence of a portfolio process $\tilde{\pi}(\cdot)$
such that \eqref{eq3} holds. According to \eqref{eu}
$$X^{x,\tilde{\pi}}(\tau)+V_{x}(\tau)=I(\mathcal{Y}_{\tau}(x)\tilde{Z}(\tau)),$$
or
\begin{equation}\label{0i9}
U'(X^{\tilde{\pi},x}(\tau)+V_{x}(\tau))=\mathcal{Y}_{\tau}(x)\tilde{Z}(\tau),
\end{equation} 
for every stopping time $\tau.$ In particular
\begin{equation}\label{pi9}
U'(X^{\tilde{\pi},x}(t)+V_{x}(t))=\mathcal{Y}_{t}(x)\tilde{Z}(t),\quad 0\leq t\leq T.
\end{equation}
 We claim that $\mathcal{Y}_{\tau}(x)=\mathcal{Y}_{0}(x)=U'(x)$ and this concludes the proof. We establish first $\mathcal{Y}_{t_{1}}(x)=\mathcal{Y}_{t_{2}}(x),$ for every $0<t_{1}\leq t_{2}\leq T. $ Arguing by contradiction assume there is $0<s<t\leq T,$ such that $\mathcal{Y}_{t}(x)\neq\mathcal{Y}_{s}(x).$ Let then $A\in\sF_s,$ and the stopping time $\hat{\tau}$ defined by
$$\hat{\tau}(\omega)\triangleq\begin{cases}t  &\text{if\,\,\,   $
\omega\in A$}\\
s   &\text{if\,\,\,   $\omega\in A^{c}$} .\end{cases}$$  
In the light of \eqref{pi9}
$$U'(X^{\hat{\pi},x}(\hat{\tau})+V(\hat{\tau}))=\mathcal{\hat{Y}}_{\hat{\tau}}(x)\tilde{Z}(\hat{\tau}),$$
where
$$\mathcal{\hat{Y}}_{\hat{\tau}}(x)\triangleq\begin{cases}\mathcal{Y}_{t}(x)&\text{if\,\,\,   $
\omega\in A$}\\
  \mathcal{Y}_{s}(x) &\text{if\,\,\,   $\omega\in A^{c}$} .\end{cases}$$
  However by \eqref{0i9} $\mathcal{\hat{Y}}_{\hat{\tau}}(x)=\mathcal{Y}_{\hat{\tau}}(x),$ a constant, hence $\mathcal{Y}_{s}(x)=\mathcal{Y}_{t}(x)=\lambda,$ a contradiction so $\mathcal{Y}_{t_{1}}(x)=\mathcal{Y}_{t_{2}}(x)=\lambda,$ for every $0< t_{1}\leq t_{2}\leq T.$ Since $t\rightarrow \mathcal{X}_{t}(\lambda)$ is continuous and $\mathcal{X}_{t}(\lambda)=x,$ then
 $\mathcal{X}_{0}(\lambda)=x,$ i.e., $\mathcal{Y}_{t}(x)=\lambda=\mathcal{Y}_{0}(x).$

  \begin{flushright}
$\diamond$
\end{flushright}

{\bf{Acknowledgements}}

\vspace{0.4cm}

The authors would like to thank Professors Tahir Choulli and Ivar Ekeland for helpful discussions and comments.

\end{document}